\documentclass[12pt]{amsart}
\usepackage{amsthm}
\usepackage{epsf}
\usepackage{amssymb}
\usepackage{latexsym}
\usepackage{amsmath}
\newtheorem{theorem}{Theorem} 
\newtheorem{definition}[theorem]{Definition}

\newtheorem{lemma}[theorem]{Lemma}

\newtheorem{proposition}[theorem]{Proposition}
\newtheorem{remark}[theorem]{Remark}

\def\Z{{\mathbb Z}}

\def\semidirect{\rtimes}

\def \ta{\tau}
\def \ta1{\tau_1}

\def\isom{{\cong}}

\newcommand\begintable[1][] {{}}

\long\def\forget#1\forgotten{}

\newif\ifXY 
\XYtrue     

\def\Z{\mathbb{Z}}

\begin{document}

\renewcommand{\subjclassname}{%
       \textup{2000} Mathematics Subject Classification}
\date{\today}

\title{Coxeter covers of the classical Coxeter groups}

\author[Amram, Shwartz, Teicher]{Meirav Amram$^{1}$, Robert Shwartz$^1$  and  Mina Teicher}
\stepcounter{footnote} \footnotetext{Partially supported by the
Emmy Noether Research Institute for Mathematics (center of the
Minerva Foundation of Germany), the Excellency Center "Group
Theoretic Methods in the Study of Algebraic Varieties" of the
Israel Science Foundation, and EAGER (EU network,
HPRN-CT-2009-00099)}.

\address{Meirav Amram, Department of Mathematics, Bar-Ilan University, Ramat-Gan 52900, Israel}
\email{meirav@macs.biu.ac.il}

\address{Robert Shwartz, Department of Mathematics, Bar-Ilan University, Ramat-Gan 52900, Israel}
\email{shwart1@macs.biu.ac.il}

\address{Mina Teicher, Department of Mathematics, Bar-Ilan University, Ramat-Gan 52900, Israel}
\email{teicher@macs.biu.ac.il}

\begin{abstract}
Let $C(T)$ be a generalized Coxeter group, which has a natural map
onto one of the classical Coxeter groups, either $B_n$ or $D_n$.
Let $C_Y(T)$ be a natural quotient of  $C(T)$, and if $C(T)$ is
simply-laced (which means all the relations between the generators
has order 2 or 3), $C_Y(T)$ is a generalized Coxeter group, too .
Let $A_{t,n}$ be a group which contains $t$ Abelian groups
generated by $n$ elements. The main result in this paper is that
$C_Y(T)$ is isomorphic to $A_{t,n} \semidirect B_n$ or $A_{t,n}
\semidirect D_n$, depends on whether the signed graph $T$
contains loops or not, or in other words C(T) is simply-laced or
not, and $t$ is the number of the cycles in $T$. This result
extends the results of Rowen, Teicher and Vishne to generalized
Coxeter groups which have a natural map onto one of the classical
Coxeter groups.
\end{abstract}


\maketitle
\section{Introduction}\label{intro}

Coxeter Groups is an  important class of groups which is used in
the study of symmetries, classifications of Lie Algebras and in
other subjects of Mathematics.

In \cite{rtv}, there is a description of Coxeter groups from which
there is a natural map onto a symmetric group. Such Coxeter groups
have natural quotient groups related to  presentations of the
symmetric group on an arbitrary set $T$ of transpositions.

These quotients, which are denoted by $C_Y(T)$, are a special type
of the generalized Coxeter groups defined in \cite{CST} by a
signed Coxeter diagram, where in addition to the regular Coxeter
relations, which arise from the graph, every signed cycle, where
the multiplication of the signs are negative, admits an extra
relation. $C_Y(T)$ is a class of groups where every negatively
signed cycle is a triangle. Hence, every extra relation has a
form: $(x_1x_2x_3x_2)^2=1$, where $x_1$, $x_2$ and $x_3$ are the
vertices of the negatively signed triangle.

The group $C_Y(T)$  also arises in the computation of certain
invariants of surfaces (see \cite{teicher}).

The paper \cite{rtv} deals with the class of Coxeter groups, whose
Coxeter diagram (the dual diagram to the diagram introduced in
\cite{rtv}) does not have  a subgraph
$\vcenter{\hbox{\epsfbox{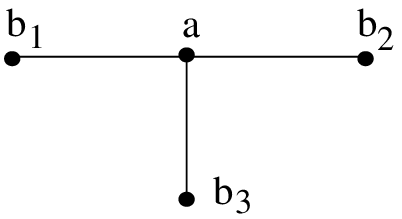}}}$ \ ($a, b_1, b_2$ and $b_3$
are Coxeter generators, $(a b_1)^3=(a b_2)^3=(a b_3)^3 = 1$ and
$(b_1 b_2)^2=(b_1 b_3)^2=(b_2 b_3)^2 = 1$) (see \cite[Remark
7.13]{rtv}).

This paper extends the results of \cite{rtv} for a wider class of
Coxeter groups $C(T)$, and $C(T)$ can be also sometimes a
generalized Coxeter group \cite{CST} where the natural
homomorphism is onto one of the classical Coxeter groups $A_n$,
$B_n$, $D_n$ (which have of course a homomorphism onto $S_n$). But
there are still Coxeter groups $C(T)$ (even simply-laced) which do
not have any homomorphism onto any of the classical Coxeter
groups, for example, $C(T)$ can not be anyone of the exceptional
Coxeter groups. In case of the configuration which mentioned above
(which is allowed in our case), two among three  vertices ($b_i$
and $b_j$) satisfy  $m_{(b_i,x)}=m_{(b_j,x)}$, for every Coxeter
generator $x$, where $m_{(b_i,x)}$ denotes the order of $b_ix$ in
$C(T)$ (the regular notation in Coxeter groups).

\vskip 0.2cm

Let us briefly recall the definitions and properties of the groups
$A_n, B_n, D_n$ and the exceptional Coxeter groups (see \cite[page
32]{hum}).  It is well known \cite{hum}, that $B_n \isom \Z_2
{\wr} S_n$ (wreath product) and $D_n$ is a subgroup of $B_n$ of
index $2$. One can present $B_n$ and  $D_n$ as groups of signed
permutations, and then present graphs of  $B_n$ and $D_n$ as
follows:
\begin{figure}[!ht]
\epsfxsize=9cm 
\epsfysize=5cm 
\begin{minipage}{\textwidth}
\begin{center}
\epsfbox{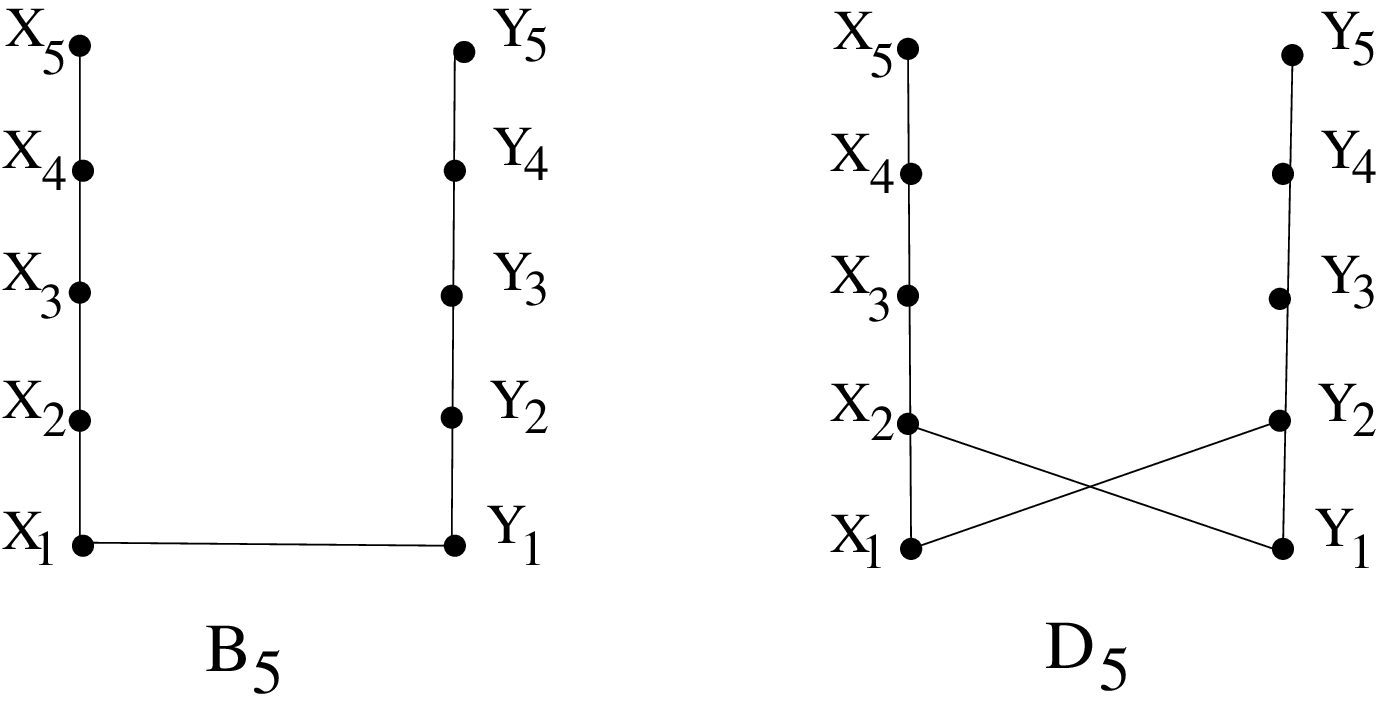}
\end{center}
\end{minipage}
\caption{}\label{dyn2}
\end{figure}
The edges in the graph which corresponds to $B_n$ are
\begin{align*}
& s_0=(x_1,y_1), s_1=(x_1,x_2)(y_1,y_2),  s_2=(x_2,x_3)(y_2,y_3),  \\
& s_3=(x_3,x_4)(y_3,y_4),  s_4=(x_4,x_5)(y_4,y_5),
\end{align*}
and the edges in the graph which correspond to $D_n$ are
\begin{align*}
& s_{\bar{1}}=(x_1,y_2)(x_2,y_1), s_1=(x_1,x_2)(y_1,y_2),
s_2=(x_2,x_3)(y_2,y_3),  \\
& s_3=(x_3,x_4)(y_3,y_4), s_4=(x_4,x_5)(y_4,y_5).
\end{align*}

We note that all the generators of $D_n$ are presented by a pair
of edges. The generators of $B_n$, apart from $s_0$, are presented
by pairs of edges, too.  This form is analogical to the $2n$
permutation presentation of  $B_n$ and  $D_n$, where $s_i$ are
presented by a product of two transpositions ($s_0$ is presented
by a single transposition in $B_n$).

\vskip 0.2cm

In Section \ref{sec2} we define the group $C(T)$ which has a
natural map onto one of the classical Coxeter groups. A diagram
for $C(T)$ (e.g. Figure \ref{ct}) is analogical to the diagram
which was introduced in \cite{rtv}, while in our case, most of
generators are presented by a couple of edges, and only specific
generators are presented by a single edge.

In Section \ref{sec3} we introduce a much more convenient
presentation of $C(T)$ by reduced diagrams. These diagrams are
signed graphs (see \cite{CST}), where the edges of the graph are
signed either by $1$ or $-1$. Signed graphs are subject to a
relation of the form ${(u_1 \cdot u_2 \cdots u_{n-1} u_n u_{n-1}
\cdots u_2)}^2=1$ for every cycle with odd number of sign $-1$
(similarly to \cite{CST}, which we call anti-cycle. Note that this
type of relations appears in \cite{CST}, but in a dual form, where
the generators are vertices and not edges. Due to this additional
relation which arises from an anti-cycle, there are signed graphs
$T$, where $C(T)$ is a generalized Coxeter group (Coxeter group
with additional relations, which arise from negatively signed
cycles, or anti-cycles). We assume that $C(T)$ is connected signed
graph, and $C(T)$ does contain a loop or at least one anti-cycle
(Otherwise the theorem is isomorphic to the Theorem in \cite{rtv}).

In Section \ref{sec4} we classify the relations which arise in the
quotient $C_Y(T)$ of $C(T)$. In addition to the anti-cyclic
relation, there are other three types of relations which arise in
$C_Y(T)$.

In Section \ref{sec5} we classify the cyclic relations, which
generate the kernel of the mapping from $C_Y(T)$ onto $B_n$ or
$D_n$. There are four possible  types of cyclic relations. Each
type defines one of the classical affine Coxeter groups,
$\tilde{A}$, $\tilde{B}$, $\tilde{C}$ and $\tilde{D}$, which are
periodic permutations or signed permutation groups (see
\cite{Eriksson}). $\tilde{A}_n$ is the well-known
$\tilde{S}_{n+1}$, where the period is $n+1$, which means
$\tilde{A}_n$ is a periodic permutation group which satisfies
$\pi(i+(n+1))=\pi(i)+(n+1)$ for every permutation in
$\tilde{A}_n$.  The other three affine Coxeter groups are periodic
sign permutations with a period of $2n+2$, which satisfies
$\pi(i+(2n+2))=\pi(i)+(2n+2)$, and in addition $\pi(-i)=-\pi(i)$,
where in the sequel, $-i$ will be denoted by $\bar{i}$, when we
treat $-1$ in a signed permutation. It is well known that
$\tilde{A}_n$ is isomorphic to $\mathbb{Z}^{n}\semidirect A_n$, or
$\mathbb{Z}^{n}\semidirect S_{n+1}$. Similarly, $\tilde{B}_n$ is
isomorphic to $\mathbb{Z}^{n}\semidirect B_n$,  $\tilde{C}_n$ is
isomorphic to $\mathbb{Z}^{n}\semidirect B_n$, $\tilde{D}_n$ is
isomorphic to $\mathbb{Z}^{n}\semidirect B_n$, where
$\mathbb{Z}^{n}$ is the group $A_{1,n}$ which will be defined in
Section \ref{sec6}.

In Section \ref{sec6} we define a group $A_{t,n}$ which will be
used for the main theorem, and in Section \ref{sec7} we prove the
main theorem which states that $C_Y(T)$ is isomorphic to the
semi-direct product of $A_{t,n}$ (which was defined in Section
\ref{sec6}) by $B_n$ or $D_n$, if the signed graph of $C(T)$
contains loops or does not contain loops, respectively.

\section{The group $C(T)$}\label{sec2}

Let $T'$ be a graph which contains $2n$ vertices $x_1, \dots, x_n$ and
$y_1, \dots, y_n$.
The edges which connect the vertices are defined as follows:
\begin{eqnarray*}
& (x_i,x_j) \ \ \ \mbox{is an edge} \iff (y_i,y_j) \ \ \ \mbox{is an edge}\\
& \mbox{and}\\
& (x_i,y_j) \ \ \ \mbox{is an edge} \iff (x_j,y_i) \ \ \ \mbox{is an edge}.
\end{eqnarray*}

\begin{figure}[!ht]
\epsfxsize=8.5cm 
\epsfysize=5cm 
\begin{minipage}{\textwidth}
\begin{center}
\epsfbox{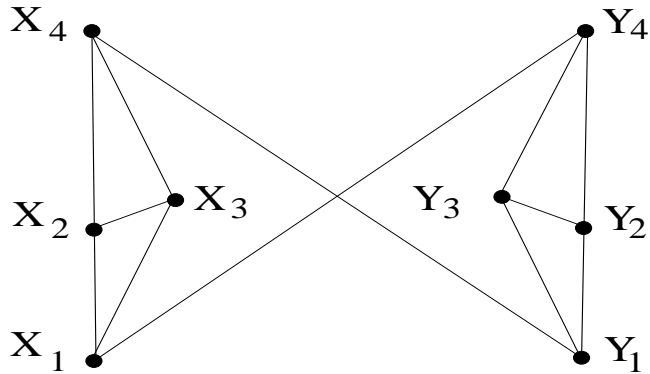}
\end{center}
\end{minipage}
\caption{An example of a graph for $C(T)$}\label{ct}
\end{figure}

For every $i \neq j$, a pair of edges $(x_i,x_j)(y_i,y_j)$ or
$(x_i,y_j)(x_j,y_i)$ presents a generator of $C(T)$. For $i = j$,
an  edge  $(x_i,y_i)$ presents a generator of $C(T)$, see for
example Figure \ref{ct}.

The group $C(T)$ admits the following relations on the edges:
\begin{enumerate}
\item[(I).] For distinct $i, j, k$ (it is the case where two pairs
of edges, which symbolize  two generators, meet at two vertices):
\begin{equation}\label{1}
((x_i,x_j)(y_i,y_j) \cdot (x_i,y_j)(x_j,y_i))^2 = 1, \ \ \
\mbox{$\vcenter {\hbox{\epsfbox{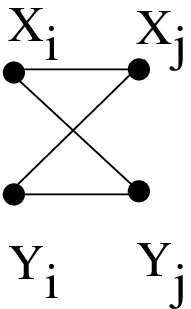}}}$}
\end{equation}
\begin{equation}\label{2}
((x_i,x_j)(y_i,y_j) \cdot (x_j,x_k)(y_j,y_k))^3 = 1, \ \ \
\mbox{$\vcenter {\hbox{\epsfbox{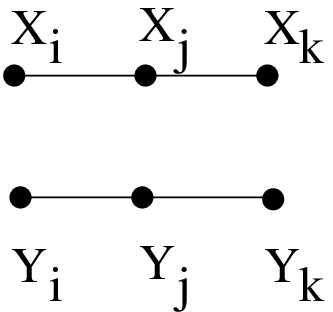}}}$}
\end{equation}
\begin{equation}\label{3}
((x_i,x_j)(y_i,y_j) \cdot (x_j,y_k)(x_k,y_j))^3 = 1,\ \ \
\mbox{$\vcenter {\hbox{\epsfbox{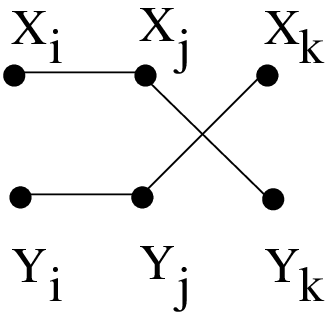}}}$}
\end{equation}
\begin{equation}\label{4}
((x_i,y_j)(x_j,y_i) \cdot (x_j,y_k)(x_k,y_j))^3 = 1, \ \ \
\mbox{$\vcenter {\hbox{\epsfbox{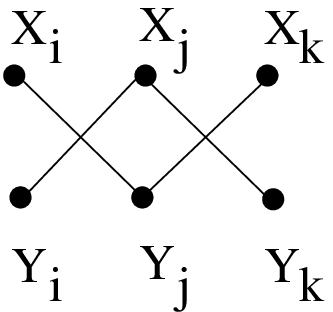}}}$}
\end{equation}
and a non simply-laced relation may hold if and only if there is a
generator of the form $(x_i,y_i)$, which admits (for distinct $i$
and $j$):
\begin{equation}\label{5}
((x_i,y_i) \cdot (x_i,x_j)(y_i,y_j))^4 = 1, \ \ \ \mbox{$\vcenter
{\hbox{\epsfbox{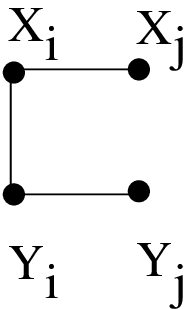}}}$}
\end{equation}
and
\begin{equation}\label{6}
((x_i,y_i) \cdot (x_i,y_j)(x_j,y_i))^4 = 1,\ \ \ \mbox{$\vcenter
{\hbox{\epsfbox{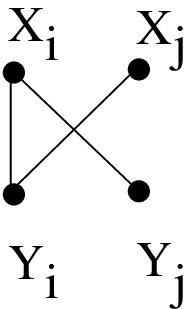}}}$}
\end{equation}
\item[(II).] For distinct $i, j, k, l$ (it is the case where two
pairs of edges are
disjoint): \\
\begin{eqnarray}
& ((x_i,x_j)(y_i,y_j) \cdot (x_k,x_l)(y_k,y_l))^2 = 1,\ \ \
\mbox{$\vcenter{\hbox{\epsfbox{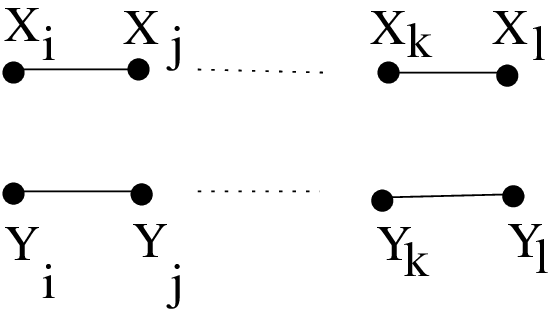}}}$} \label{7} \\
& ((x_i,y_j)(x_j,y_i) \cdot (x_k,x_l)(y_k,y_l))^2 = 1,\ \ \
\mbox{$\vcenter{\hbox{\epsfbox{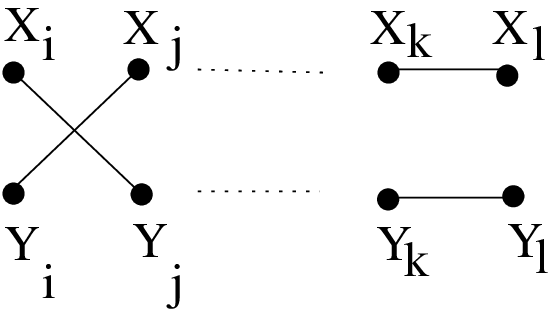}}}$} \label{8} \\
& ((x_i,y_j)(x_j,y_i) \cdot (x_k,y_l)(x_l,y_k))^2 = 1, \ \ \
\mbox{$\vcenter{\hbox{\epsfbox{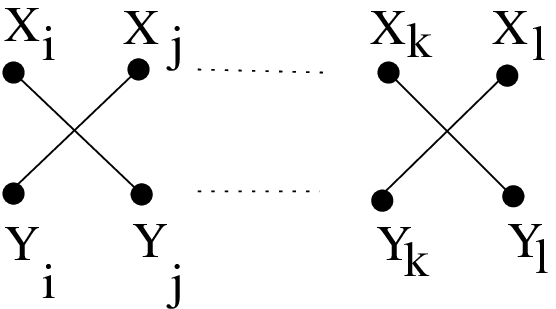}}}$} \label{9}
\end{eqnarray}
\item[(III).] For distinct $i$, $j$ and $k$  (it is the case where
an edge $(x_i,y_i)$
is disjoint from a pair of edges):\\
\begin{eqnarray}
& ((x_i,y_i) \cdot (x_j,x_k)(y_j,y_k))^2 = 1, \ \ \
\mbox{$\vcenter{\hbox{\epsfbox{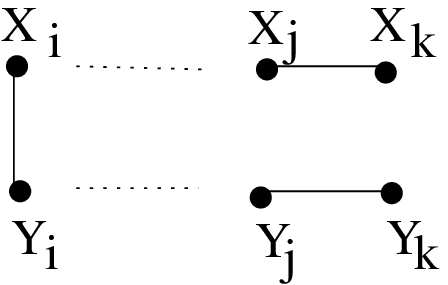}}}$} \label{10} \\
& ((x_i,y_i) \cdot (x_j,y_k)(x_k,y_j))^2 = 1, \ \ \
\mbox{$\vcenter{\hbox{\epsfbox{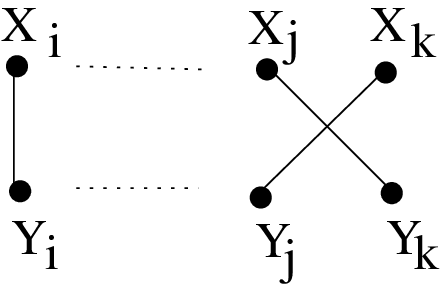}}}$} \label{11} \\
& ((x_i,y_i) \cdot (x_j,y_j))^2 = 1, \ \ \
\mbox{$\vcenter{\hbox{\epsfbox{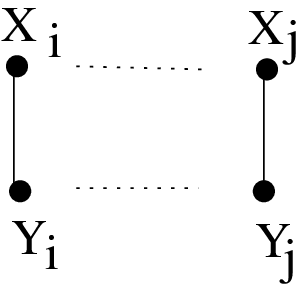}}}$} \label{12}
\end{eqnarray}
\end{enumerate}

Each graph $T'$ which satisfies the above described relations, has a
natural mapping into $B_n$ or $D_n$. Each pair of the form
$(x_i,x_j)(y_i,y_j)$ is mapped to the element
$(ij)(\bar{i}\bar{j})$; a pair of the form $(x_i,y_j)(x_j,y_i)$ is
mapped to the element  $(i \bar{j})(\bar{i} {j})$; and an edge of
the form $(x_i,y_i)$  is mapped to the transposition  $(i
\bar{i})$. In the case that there are no edges of the form
$(x_i,y_i)$, the group $C(T)$ has a natural map into $D_n$ (and
$C(T)$ is simply-laced).

\section{The reduced signed graphs}\label{sec3}

Due to the symmetry between $x_i$ and $y_i$, we may consider an equivalent
reduced signed graph $T$.

Instead of a graph $T'$ with $2n$ vertices, we consider a signed graph $T$
\cite{CST} with only $n$ vertices, such that there are two types
of edges, which connect the vertices. We replace
$(x_i,x_j)(y_i,y_j)$ by $(x_i,x_j)_1$, and  $(x_i,y_j)(x_j,y_i)$
by   $(x_i,x_j)_{-1}$. We replace also  $(x_i,y_i)$ by a loop
$(x_i,x_i)_{-1}$.

Then  $B_n$ and $D_n$ are presented by graphs in Figure
\ref{dyn15} (see original graphs in Figure \ref{dyn2} for
comparison):
\begin{figure}[!ht]
\epsfxsize=8.5cm 
\epsfysize=5cm 
\begin{minipage}{\textwidth}
\begin{center}
\epsfbox{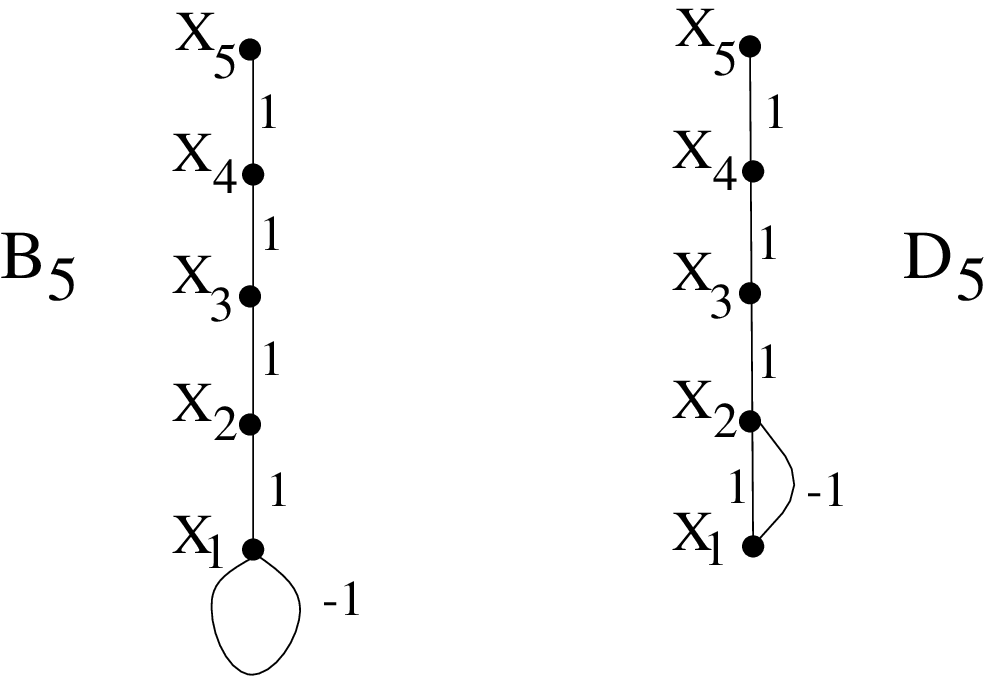}
\end{center}
\end{minipage}
\caption{}\label{dyn15}
\end{figure}

We note that the type of the group $C(T)$ in \cite{rtv} can be
presented as a graph, where all edges are of type $1$. This is due
to the existence of a natural mapping of $C(T)$ onto the symmetric
group $S_n$.

\begin{definition}
The edges  $(x_i,x_j)_1$ and  $(x_i,x_j)_{-1}$ are called conjugated
edges.
\end{definition}

The relations which hold in the reduced graph are induced from the
ones, which relate to the original graphs. For $a, b \in \{1,
-1\}$ and for distinct $i, j, k, l$ we have:\\
two conjugated edges commute
\begin{equation}\label{13}
((x_i,x_j)_1 \cdot (x_i,x_j)_{-1})^2 = 1 \ \ \mbox{derived from
(\ref{1})},
\end{equation}
two edges meet at a vertex
\begin{equation}
((x_i,x_j)_a \cdot (x_j,x_k)_b)^3 = 1 \ \ \mbox{derived from
(\ref{2})-(\ref{4})},
\end{equation}
a loop and an edge meet at a vertex
\begin{equation}
((x_i,x_i)_{-1} \cdot (x_i,x_j)_a)^4 = 1 \ \ \mbox{derived from
(\ref{5})-(\ref{6})},
\end{equation}
two edges are disjoint
\begin{equation}
((x_i,x_j)_a \cdot (x_k,x_l)_b)^2 = 1 \ \ \mbox{derived from
(\ref{7})-(\ref{9})},
\end{equation}
a loop and an edge are disjoint
\begin{equation}
((x_i,x_i)_{-1} \cdot (x_j,x_k)_a)^2 = 1 \ \ \mbox{derived from
(\ref{10})-(\ref{11})},
\end{equation}
two loops are disjoint
\begin{equation}
((x_i,x_i)_{-1} \cdot (x_j,x_j)_{-1})^2 = 1 \ \ \mbox{derived from
(\ref{12})}.
\end{equation}

In addition there is a relation which arises from cycles with odd
number of edges, signed by $-1$, similarly to the relation which
appears in \cite[Page 193]{CST}. We call it an anti-cycle
relation.

\begin{definition}{\bf{Anti-cycles}}\\
Let $x_1, \dots, x_n$ be $n$ vertices on a cycle. The edges are
\begin{equation*}
u_1:=(x_1,x_2)_{a_1}, \ \ \dots, \ \
u_{n-1}:=(x_{n-1},x_n)_{a_{n-1}}, \ \ u_n:=(x_n,x_1)_{a_n},
\end{equation*}
where $a_i \in \{1, -1\}$, $1 \leq i \leq n$ and $\#\{a_i \ | \
a_i=-1\} $ is odd.
\end{definition}
In this case we have:
\begin{equation}
(u_1 u_2 \cdots u_{n-1} \cdot u_n u_{n-1} \cdots u_2)^2=1.
\end{equation}

In a similar way, we derive relations of the form (for $1 \leq i
\leq n$)
\begin{equation}
(u_i \cdot  u_{i+1} \cdots u_n u_1 \cdots u_{i-1} u_{i-2} u_{i-3}
\cdots u_1 u_n \cdots u_{i+1})^2=1.
\end{equation}

\begin{remark}\label{anti}
If a signed graph $T$ does not contain any anti-cycle (even no
conjugated edges, which is an anti-cycle of length two) neither a loop,
then the graph $T$ describes the same groups which appears in
\cite{rtv}, where the natural homomorphism is by omiting the signs. It is
homomorphism, since the additional relation which described in this paper
caused by anti-cycle relations (including conjugated edges) or by
relations involving loops. Hence, we assume that $T$ contains at least
one anti-cycle or a loop (otherwise the result is in \cite{rtv}).
\end{remark}

\vskip 0.2cm

\begin{figure}[!ht]
\epsfxsize=5cm 
\epsfysize=4cm 
\begin{minipage}{\textwidth}
\begin{center}
\epsfbox{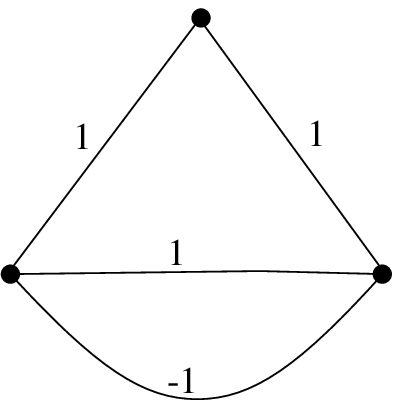}
\end{center}
\end{minipage}
\caption{}\label{clown}
\end{figure}

There are graphs $T$ where this additional relation makes $C(T)$
to be a generalized Coxeter Group as it appears in \cite{CST}. For
example, in Figure \ref{clown} one can find a group, which is a
generalized one, since we have an anti-cycle and a cycle, which
contain the same three vertices.

\vskip 0.2 cm

\begin{remark}
 We notice that the most simple case for an
anti-cycle are two conjugated edges $u_1=(x_1,x_2)_1$ and
$u_2=(x_1,x_2)_{-1}$ which form an anti-cycle. Then the relation
is just saying $u_1$ commutes with $u_2$ which we have already
assumed (see Relation (\ref{13})).
\end{remark}

\begin{lemma}\label{natural}
Let $T$ be a connected signed graph with $n$ vertices $x_1, \cdots
x_n$, and let $\phi: C(T)\rightarrow B_n$ the natural mapping such that
$\phi((x_i x_j)_1)=(i j)(\bar{i} \bar{j})$,
$\phi((x_i x_j)_{-1})=(\bar{i} j)(i \bar{j})$, and $\phi((x_i
x_i)_{-1})=(i \bar{i})$ for every $1\leq i,j\leq n$. Then the following
holds:

1) If $T$ does not contain a loop nor an anti-cycle then $Im(\phi)$ is
a subgroup of $B_n$ isomorphic to $S_n$.

2) If $T$ does contain an anti-cycle but does not contain a loop, then
$Im(\phi)=D_n$.

3) If $T$ does contain a loop then $Im(\phi)=B_n$.
\end{lemma}

We use three propositions to prove the lemma.
\begin{proposition}\label{conj}
Let $x_1\cdots x_k$ be $k$ vertices in an anti-cycle, where the edges are
$w_i:=(x_{i-1} x_i)_{a_{i-1}}$ and $w_1:=(x_k x_1)_{a_k}$. Then

$\phi(w_{i+1}w_{i+2}\cdots w_{k}w_{1}\cdots w_{i-2}w_{i-1}w_{i-2}\cdots
w_{1}w_{k}\cdots w_{i+1})=(i-1, i)(\overline{i-1}, \bar{i})$, in
case $a_{i-1}={-1}$ which means, $\phi(w_i)=(\overline{i-1}, i)(i-1,
\bar{i})$.

$\phi(w_{i+1}w_{i+2}\cdots w_{k}w_{1}\cdots w_{i-2}w_{i-1}w_{i-2}\cdots
w_{1}w_{k}\cdots w_{i+1})=(\overline{i-1}, i)(i-1, \bar{i})$
in case $a_{i-1}=1$ which means, $\phi(w_i)=(i-1, i)(\overline{i-1},
\bar{i})$.
\end{proposition}

\begin{proposition}\label{conje}
Let be a signed path connected to an anti-cycle, where
$x_1\cdots x_k$ be $k$ vertices in an anti-cycle, and the edges are
$w_i:=(x_{i-1} x_i)_{a_{i-1}}$ and $w_1:=(x_k x_1)_{a_k}$ and the vertices
of the path are $x_k,\cdots, x_s$ and the connecting edges are
$w_i:=(x_{i-1} x_{i})_{a_i}$ for $k+1\leq i\leq s$. Then

$\phi(w_i^{w_{i-1}\cdots w_{k+1}w_k\cdots w_2\cdots w_kw_1w_{k+1}\cdots
w_{i-1}})=(i-1, i)(\overline{i-1}, \bar{i})$ in
case $a_{i-1}=-1$ which means, $\phi(w_i)=(\overline{i-1}, i)(i-1,
\bar{i})$. and

$\phi(w_i^{w_{i-1}\cdots w_{k+1}w_k\cdots w_2\cdots w_kw_1w_{k+1}\cdots
w_{i-1}})=(\overline{i-1}, i)(i-1, \bar{i})$ in
case $a_{i-1}=1$ which means, $\phi(w_i)=(i-1, i)(\overline{i-1},
\bar{i})$.

where $a^b$ means $a$ conjugated by $b$.
\end{proposition}

\begin{proposition}\label{conjl}
Let be a signed path connected to a loop, where $x_0$ is a vertex
containing a loop $v$, and $w_i:=(x_{i-1} x_i)_{a_{i-1}}$ are the vertices
of a path. Then

$\phi(w_{i-1}\cdots w_1vw_1\cdots w_{i-1}w_iw_{i-1}\cdots w_1vw_1\cdots
w_{i-1})=(i-1, i)(\overline{i-1}, \bar{i})$ in
case $a_{i-1}=-1$ which means, $\phi(w_i)=(\overline{i-1}, i)(i-1,
\bar{i})$.

$\phi(w_{i-1}\cdots w_1vw_1\cdots w_{i-1}w_iw_{i-1}\cdots w_1vw_1\cdots
w_{i-1})=(i-1, i)(\overline{i-1}, \bar{i})$ in
case $a_{i-1}=1$ which means, $\phi(w_i)=(i-1, i)(\overline{i-1}, \bar{i})$.
\end{proposition}

\vskip 0.2cm \noindent {\bf{Proof of Lemma \ref{natural}}}
Assume 1) holds. Then $T$ does not contain a loop nor an anti-cycle,
then by omiting the signs of $T$, mapping the edges onto $S_n$ (remark
\ref{anti}).

Assume 2) holds. Since $T$ is connected and contains at least one
anti-cycle, every edge in $T$ either lies on an anti-cycle or connected
by a path to an anti-cycle. Hence, if $\phi((x_i x_j)_1)=(i j)(\bar{i}
\bar{j})$, then by Propostions \ref{conj} and \ref{conje} there exists an
element $w\in C(T)$ such that $\phi(w)=(\bar{i} j)(i \bar{j})$. On the
other hand, if $\phi((x_i x_j)_{-1})=(\bar{i} j)(i \bar{j})$ then by the
same argument there exists $w$ such that $\phi(w)=(i j)(\bar{i} \bar{j})$.
Since $T$ is connected, there is a path connecting any two vertices in
$T$, then by the same argument as in \cite{rtv} for every
distinct $i$ and $j$ such that $1\leq i,j\leq n$, there are elements
$w_1$ and $w_2$ such that $\phi(w_1)=(i j)(\bar{i} \bar{j})$ and
$\phi(w_2)=(\bar{i} j)(i \bar{j})$. The subgroup of $B_n$ which is
generated by all signed transpositions is $D_n$.

Assume 3 holds. Since $T$ is connected and contains a loop, every edge
which is not a loop connected with a path to a loop.
Hence, if $\phi((x_i x_j)_1)=(i j)(\bar{i}
\bar{j})$, then by Propostion \ref{conjl} there exists an
element $w\in C(T)$ such that $\phi(w)=(\bar{i} j)(i \bar{j})$. On the
other hand, if $\phi((x_i x_j)_{-1})=(\bar{i} j)(i \bar{j})$ then by the
same argument there exists $w$ such that $\phi(w)=(i j)(\bar{i} \bar{j})$.
Since $T$ contains a loop, then there exists an element $v$ such that
$\phi(v)=(i \bar{i})$, and the subgroup of $B_n$ which is generated by
allthe signed transpositions $(i j)(\bar{i} \bar{j})$, $(\bar{i} j)(i
\bar{j})$ and an element of a form $(i \bar{i})$ is all $B_n$.

\section{The group $C_Y(T)$}\label{sec4}
We define the group $C_Y(T)$ as a quotient of $C(T)$ by the 'fork'
relations. The fork relations in $C(T)$ are (for $a, b, c
\in \{1, -1\}$):  \\
I. \ Three edges meet at a common vertex: \ \ \ \ \
$\vcenter{\hbox{\epsfbox{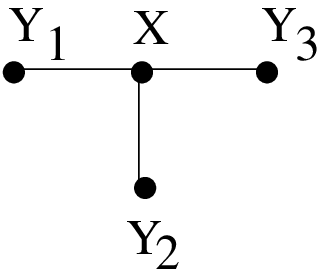}}}$.
$$
((x,y_1)_a \cdot (x,y_2)_b)^3 = ((x,y_1)_a \cdot (x,y_3)_c)^3 = ((x,y_2)_b
\cdot (x,y_3)_c)^3 = 1.
$$
Then ($R_1$) is  (as in \cite{rtv}):
\begin{equation}\label{r1}
((x,y_1)_a \cdot (x,y_2)_b (x,y_3)_c  (x,y_2)_b)^2 = 1.
\end{equation}
II. \ Two conjugated edges $(x_2,x_3)_1$ and $(x_2,x_3)_{-1}$ meet
at both of their common vertices ($x_2$ and $x_3$), two other
edges $(x_1,x_2)_a$ and $(x_3,x_4)_b$ \ \ \
$\vcenter{\hbox{\epsfbox{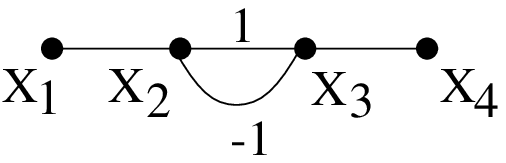}}}$. \\
Then ($R_2$)  is:
\begin{equation}\label{r2}
((x_1,x_2)_a  (x_2,x_3)_1 (x_1,x_2)_a \cdot (x_3,x_4)_b
(x_2,x_3)_{-1} (x_3,x_4)_b)^2 = 1.
\end{equation}
III. \ A loop and two edges meet at a vertex \ \ \
$\vcenter{\hbox{\epsfbox{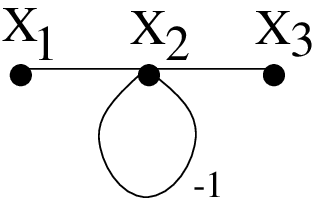}}}$. \\
Then ($R_3$)  is:
\begin{equation}{(R_3)}\label{r3}
((x_2,x_2)_{-1} \cdot (x_1,x_2)_a (x_2,x_3)_b (x_1,x_2)_a)^2 =1,
\end{equation}
and ($R_4$)  is:
\begin{equation}\label{r4}
((x_1,x_2)_a \cdot (x_2,x_2)_{-1} (x_2,x_3)_b (x_2,x_2)_{-1})^3
=1.
\end{equation}
\medskip

\vskip 0.2cm

We recall that in order to prove these relations, we  consider
$u_i$ as a signed permutation in $B_n$, where $(x_i,x_{i+1})_1$ is
$(i, {i+1})(\bar{i}, \overline{i+1})$ and $(x_i,x_{i+1})_{-1}$ is
$(i, \overline{i+1})(\bar{i}, {i+1})$.

Note that in the case of $D$-covers, we may have only (\ref{r1})
and (\ref{r2}), since (\ref{r3}) and (\ref{r4}) involve loops,
which may appear only in $B$-covers. Thus:
$$
C_Y(T) = C(T) / \langle{(\ref{r1})  \cup (\ref{r2}) }\rangle \ \ \
\mbox{for D-covers}
$$
and
$$
C_Y(T) = C(T) / \langle{(\ref{r1})  \cup (\ref{r2})  \cup
(\ref{r3}) \cup (\ref{r4}) } \rangle\ \ \ \mbox{for B-covers}.
$$

\section{Mapping $C_Y(T)$ onto $B_n$ or $D_n$}\label{sec5}
Now we classify the relations, which may appear in the kernel of
the mapping from $C_Y(T)$ onto $B_n$ or $D_n$ (similarly as done
for the 'cyclic' relations in \cite{rtv}).

\begin{enumerate}
\item[(I).] Cycles: \nonumber \\
Let $T$ be connected signed graph which contains at least one
anti-cycle. Let $x_0, \dots, x_{m-1}$ be $m$ vertices on a cycle, which
are connected by the $m$ edges $(x_{i-1},x_i)_{a_{i-1}}$, and
$(x_{m-1},x_0)_{a_{m-1}}$ where
$\#\{a_i \ | \ a_i=-1\} $ is even.

If $a_{i-1}=1$, then $u_i:=(x_{i-1} x_i)_{a_{i-1}}$.

If $a_{i-1}=-1$, then $\bar{u}_{i}:=(x_{i-1} x_i)_{a_{i-1}}$.

Now define $u_i$ for the cases where $a_{i-1}=-1$, and $\bar{u}_{i}$ for
the cases where $a_{i-1}=1$. Since, $T$ is connected and $T$ does
contain an anti-cycle or a loop, let $w_1,\cdots w_k$ be $k$ edges which
form an anti-cycle of length $k$ in case $T$ contains an anti-cycle,
otherwise, lew $w$ be a loop. Let $v_1, \cdots v_s$ be a path connecting
the anti-cycle of length $k$ or the loop with the cycle of length $m$.
Then:

\vskip 0.2cm

In case $a_0=-1$:

\begin{equation}\label{cyc1}
u_1:= \bar{u}_{1}^{v_{s}\cdots
v_{1}w_{1}w_k^{w_{k-1}\cdots w_2}v_{1}\cdots v_{s}}
\end{equation}

and in case $a_0=1$:

\begin{equation}\label{cycb1}
\bar{u}_1:= u_{1}^{v_{s}\cdots
v_{1}w_{1}w_k^{w_{k-1}\cdots w_2}v_{1}\cdots v_{s}}
\end{equation}

Then inductively we define $u_i$ where $a_{i-1}=-1$ and $\bar{u}_i$ where
$a_{i-1}=1$ for every $1\leq i\leq m$ as following

\begin{equation}\label{cyci}
u_i:= \bar{u}_{i}^{u_{i-1}\cdots u_{1}v_{s}\cdots
v_{1}w_{1}w_k^{w_{k-1}\cdots w_2}v_{1}\cdots v_{s}u_{1}\cdots u_{i-1}}.
\end{equation}

\begin{equation}\label{cycbi}
\bar{u}_i:= u_{i}^{u_{i-1}\cdots u_{1}v_{s}\cdots
v_{1}w_{1}w_k^{w_{k-1}\cdots w_2}v_{1}\cdots v_{s}u_{1}\cdots u_{i-1}}
\end{equation}

where we denote $a^b$ instead of $b^{-1}ab$.

\vskip 0.2cm

In case of loop instead of anti-cycle, we write $w$ instead of
$w_1w_k^{w_{k-1}\cdots w_2}$ in equations \ref{cyc1}, \ref{cycb1},
\ref{cyci} and \ref{cycbi}.

\begin{remark}
We notice that the existence of an anti-cycle or a loop connecting the
cycle allows us to define $u_i$ and $\bar{u}_i$ for every $1\leq i\leq n$,
such that the natural mapping $\phi$ from $C(T)$ onto $B_n$ or $D_n$
satisfies

$\phi(u_i)=(i-1 ,i)(\overline{i-1}, \bar{i})$ and
$\phi(\bar{u}_{i})=(i-1,
\bar{i})(\overline{i-1}, i)$, $\phi(u_m)=(m-1, 0)(\overline{m-1},
\bar{0})$ and $\phi(\bar{u}_{m})=(m-1, \bar{0})(\overline{m-1}, 0)$.
\end{remark}

\medskip
\item[(II).] $\tilde{D}$-type cycles: \nonumber \\
Two anti-cycles connected by a path are called a $\tilde{D}$-type
cycle. The length of the anti-cycles can be every length $\geq 2$
(anti-cycle of length 2 means two conjugated edges). Let $x_0,
\dots, x_{m-1}$ be $m$ vertices, where $x_0, \dots, x_{k_1-1}$
form an anti-cycle of length $k_1$, and $x_{k_2}, \dots, x_{m-1}$
form another anti-cycle, and there is a simple signed path
connecting the vertices $x_{k_1-1}$ and $x_{k_2}$.  We define
$u_i$ and $\bar{u}_{i}$ for every $1\leq i\leq m-1$. In case
$k_1=2$: $u_1:=(x_0,x_1)_{1}$, $\bar{u}_{1}:=(x_0,x_1)_{-1}$,
otherwise we look at the sign of the edge connecting $x_0$ and
$x_1$. If the sign is $+1$ then $u_1:={(x_0,x_1)}_{1}$ and
$\bar{u}_{1}:={(x_{k_1-1},x_0)}_{a_{k_1-1}}$ conjugated by
${(x_1,x_2)}_{a_1}{(x_2,x_3)}_{a_2} \cdots
{(x_{k_1-2},x_{k_1-1})}_{a_{k_1-2}}$, where $a_i$ is the sign of
the edge connecting $x_i$ with $x_{i+1}$. If the sign of the edge
connecting $x_0$ with $x_1$ is $-1$, then
$\bar{u}_{1}:={(x_0,x_1)}_{-1}$ and
$u_1:={(x_{k_1-1},x_0)}_{a_{k_1-1}}$ conjugated by
${(x_1,x_2)}_{a_1}{(x_2,x_3)}_{a_2}\cdots
{(x_{k_1-2},x_{k_1-1})}_{a_{k_1-2}}$.

Similarly, we define $u_{m-1}$ and $\bar{u}_{m-1}$ where we look
at the second anti-cycle. If the length of the second anti-cycle
is $2$, then:
$$u_{m-1}:={(x_{m-2},x_{m-1})}_{1}, \ \ \bar{u}_{m-1}:={(x_{m-2},x_{m-1})}_{-1},$$
otherwise, similarly to the definition of $u_1$ and $\bar{u}_1$ we
look at the sign of the edge connecting $x_{m-2}$ and $x_{m-1}$.
If the sign is $+1$ then $u_{m-1}:={(x_{m-2},x_{m-1})}_{1}$ and
$\bar{u}_{m-1}:={(x_{k_2},x_{m-1})}_{a_{m-1}}$ conjugated by
${(x_{m-1},x_{m-2})}_{a_{m-2}}{(x_{m-2},x_{m-3})}_{a_{m-3}}\cdots
{(x_{k_2+1},x_{k_2})}_{a_{k_2}}$, where $a_i$ is the sign of the
edge connecting $x_i$ with $x_{i+1}$. If the sign of the edge
connecting $x_0$ with $x_1$ is $-1$, then
$\bar{u}_{m-1}:={(x_{m-2},x_{m-1})}_{-1}$ and
$u_{m-1}:={(x_{k_2},x_{m-1})}_{a_{m-1}}$ conjugated by
$(x_{m-1},x_{m-2})_{a_{m-2}}{(x_{m-2},x_{m-3})}_{a_{m-3}}\cdots
{(x_{k_2+1},x_{k_2})}_{a_{k_2}}$.

And  we define $u_i$ and $\bar{u}_i$ for every $2\leq i\leq m-2$
in the following way. We denote an edge in the signed graph (for
$2 \leq i \leq m-2$) as ${(x_{i-1},x_i)}_{a_i}$.

If $a_i=1$, then:
$$u_i:={(x_{i-1},x_i)}_1$$
and
$$\bar{u}_i:=u_{i-1} u_{i-2} \cdots u_2 u_1 \bar{u}_1 u_2
\cdots u_{i-1} u_i u_{i-1} \cdots u_2 u_1 \bar{u}_1 u_2 \cdots
u_{i-2} u_{i-1}.$$
If $a_i=-1$, then:
$$\bar{u}_i:={(x_{i-1},x_i)}_{-1},$$
and
$$u_i:=u_{i-1}u_{i-2}\cdots u_2 u_1 \bar{u}_1  u_2 \cdots
u_{i-1}\bar{u}_iu_{i-1}\cdots u_2 u_1 \bar{u}_1u_2\cdots
u_{i-2}u_{i-1}.$$
\begin{figure}[!ht]
\epsfxsize=9cm 
\epsfysize=4cm 
\begin{minipage}{\textwidth}
\begin{center}
\epsfbox{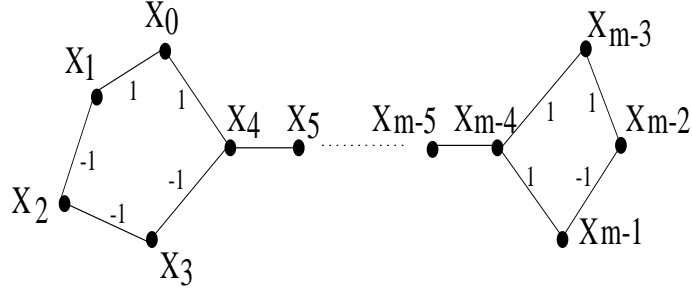}
\end{center}
\end{minipage}
\caption{$\tilde{D}$-type cycle}\label{star1}
\end{figure}

Moreover, we define  elements $u_m$ and $\bar{u}_m$ to be
$$u_m:= \bar{u}_{1}
u_{2} u_{3} \cdots u_{m-2} \bar{u}_{m-1} u_{m-2} \cdots u_{3}
u_{2} \bar{u}_{1}$$ and
$$\bar{u}_m:= u_{1} u_{2} u_{3} \cdots
u_{m-2} \bar{u}_{m-1} u_{m-2} \cdots u_{3} u_{2} u_{1}.$$

\medskip
\item[(III).] $\tilde{B}$-type cycles: \nonumber \\
A loop and an  anti-cycle which are connected by a path are called
$\tilde{B}$-type cycle. The length of the anti-cycles can be every
length $\geq 2$.

Let $x_0, \dots, x_{m-1}$ be $m$ vertices, where we have a loop in
$x_0$, an anti-cycle connecting the vertices $x_k$ and $x_{m-1}$,
and a simple signed path between $x_0$ and $x_k$.

We define $u_i$ and $\bar{u}_i$ in the following way (for $1\leq
i\leq m-1$):

Let $v:=(x_0,x_0)_{-1}$. If $(x_0,x_1)_1$ belongs to the signed
graph, then $u_1:=(x_0,x_1)_1$ and $\bar{u}_1:=vu_1v$. Otherwise,
for $(x_0,x_1)_{-1}$ belonging to the signed graph,
$\bar{u}_1:=(x_0,x_1)_{-1}$ and $u_1:=v\bar{u}_1v$.

For $2\leq i\leq m-1$, we define $u_i$ in the same way as it was
defined for $\tilde{D}$-type cycles.
\begin{figure}[!ht]
\epsfxsize=9cm 
\epsfysize=4cm 
\begin{minipage}{\textwidth}
\begin{center}
\epsfbox{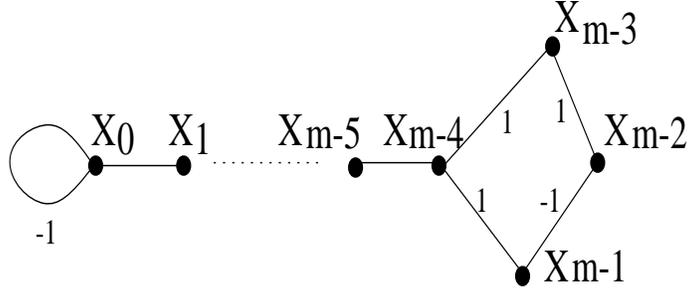}
\end{center}
\end{minipage}
\caption{$\tilde{B}$-type cycle}\label{star2}
\end{figure}

 Moreover, we define elements $u_m$ and $\bar{u}_m$ as follows:
$$u_m:= v u_{1} u_{2} \cdots u_{m-2} \bar{u}_{m-1} u_{m-2} \cdots u_{2}
u_{1} v$$
and
$$\bar{u}_m:= u_{1} u_{2} \cdots u_{m-2} \bar{u}_{m-1} u_{m-2} \cdots u_{2} u_{1}.$$

\begin{proposition}\label{mapping}
Consider the natural mapping $\phi$ from the
$\tilde{D}$-type  or $\tilde{B}$-type cycle of length $m$ onto $D_m$,
or $B_m$ $\phi(u_i)=(i-1 ,i)(\overline{i-1}, \bar{i})$ and
$\phi(\bar{u}_{i})=(i-1,
\bar{i})(\overline{i-1}, i)$, $\phi(u_m)=(m-1, 0)(\overline{m-1},
\bar{0})$ and $\phi(\bar{u}_{m})=(m-1, \bar{0})(\overline{m-1}, 0)$
\end{proposition}

\begin{remark}
We notice that there is an edge on $\tilde{D}$-type or
on $\tilde{B}$-type cycle which one $u_{m-1}$ admits only from the
defined $u_i$ and $\bar{u}_{i}$ (The edge connecting $x_{m-2}$ to
$x_{m-1}$ or the edge connecting $x_{m-1}$ to $x_{m-4}$ depends on
$a_{m-2}$). This edge will be important when we define the spanning `tree`
in section \ref{sec7}, where we omit this
edge. Hence, by omiting this edge we omit $u_{m-1}$ only.

By symmetry we can define $u_i$ in diferent way, such there will be one
edge only in one of the anti-cycles such that one of the $u_i$'s admits
only, and $\phi(u_i)$ satisfies the condinitions of Proposition
\ref{mapping}.
\end{remark}

\medskip
\item[(IV).] $\tilde{C}$-type cycles: \nonumber \\
Two loops connected by a simple path are called a $\tilde{C}$-type
cycle. Let $x_0, \dots, x_{m-1}$ be $m$ vertices, and two loops
$$v:=(x_0,x_0)_{-1}, \ \ \ \ w:=(x_{m-1},x_{m-1})_{-1}.$$
We define $u_i$ in the same way as it was defined for
$\tilde{B}$-type cycles ($1\leq i\leq m-1$).  \ \ \
\begin{figure}[!ht]
\epsfxsize=9cm 
\epsfysize=2cm 
\begin{minipage}{\textwidth}
\begin{center}
\epsfbox{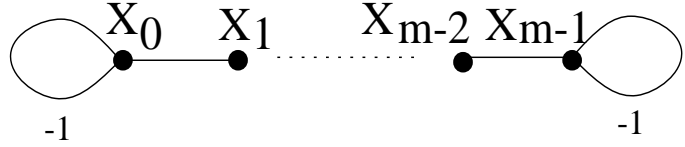}
\end{center}
\end{minipage}
\caption{$\tilde{C}$-type cycle}\label{dyn21}
\end{figure}

\end{enumerate}

In addition we define elements $u_m$ and $\bar{u}_m$ in the
following way:
$$u_m:= u_{1} u_{2} \cdots u_{m-2} u_{m-1} u_{m-2}
\cdots u_{2} u_{1} $$ and $$\bar{u}_m:= u_{1} u_{2} \cdots
u_{m-2} w u_{m-1} w u_{m-2} \cdots u_{2} u_{1}.$$

\begin{remark}
Propositon \ref{mapping} holds for the natural mapping from
$\tilde{C}$-type cycle of length $m$ onto $B_m$ too. We notice that
from the defined elements  $u_i$ and $\bar{u}_{i}$, the element
$\bar{u_m}$ admits only the loop $w$, where in the spanning `tree` (will
be defined in section \ref{sec7} we
omit this loop, hence omiting again $\bar{u}_{m}$ only from the defined
$u_i$'s and $\bar{u}_{i}$'s
\end{remark}

\begin{proposition}
Let $\phi$ be the natural map from one of the cycle onto $B_n$ or $D_n$.
Then:
$$\phi(u_1 u_2 u_3 \cdots u_{m-1}) = \phi(u_2 u_3 u_4 \cdots u_{m})$$
and
$$\phi(u_{m}u_{m-1}\bar{u}_{m-1}u_{m}u_1 u_2 u_3 \cdots
u_{m-3}\bar{u}_{m-2}u_{m-1}) = \phi(u_1u_{m}\bar{u}_{m}u_1u_2 u_3 u_4
\cdots u_{m-2}\bar{u}_{m-1}u_m).$$
\end{proposition}

\begin{proof}
The first equation has been proved in \cite{rtv}, and the second
one we get easily by substituting the signed permutation
$\phi(u_i)$ where $\phi$ is the natural map from $C_Y(T)$ onto
$B_n$ or $D_n$. By Proposition \ref{mapping} for $1\leq i\leq m-1$,
$\phi(u_i)=(i-1,i)(\overline{i-1},\bar{i})$ and
$\phi(\bar{u_i})=(\overline{i-1},i)(i-1,\bar{i})$, and
$\phi(u_m)=(m-1,0)(\overline{m-1},\bar{0})$,
$\phi(\bar{u_m}=(\overline{m-1},0)(m-1,\bar{0})$. Then:
\begin{eqnarray*}
\phi(u_{m}u_{m-1}\bar{u}_{m-1}u_{m}u_1 u_2 u_3 \cdots
u_{m-3}\bar{u}_{m-2}u_{m-1}) = \phi(u_1u_{m}\bar{u}_{m}u_1u_2 u_3
u_4 \cdots u_{m-2}\bar{u}_{m-1}u_m)\\
=(m-1 \; \;  m-2 \dots 0)(\overline{m-1} \; \; \overline{m-2}
\dots \bar{0}).
\end{eqnarray*}
\end{proof}

The definition of $u_i$ and $\bar{u}_i$ for $1\leq i\leq m$
are important, since it enables defining $\gamma_{i}$ and
$\gamma_{\bar{i}}$ for every $1\leq i\leq m$ for every type
($\tilde{A}$ or $\tilde{B}$ or $\tilde{C}$ or $\tilde{D}$) of
cycle which contains $n$ vertices, as defined in \cite{rtv} (see
Section \ref{sec5}).

We call  the above figures $\tilde{B}$-, $\tilde{C}$- and
$\tilde{D}$-types cycles, since the groups which are described by
them are the affine groups $\tilde{B}$, $\tilde{C}$ and
$\tilde{D}$. By \cite{MV}, an infinite Coxeter group is large if
and only if the group is not affine. Hence, a diagram $T$ defines a
large group $C_Y(T)$ (quotient of $C(T)$ by one of the relations which
are mentioned in Section \ref{4}) for every graph other than one of the
cycles which are mentioned here,  an anti-cycle (which is a graph of
$D_n$), a line connecting an anti-cycle or a loop. In Section \ref{7} we
will conclude that $C_Y(T)$ is large if and only if $T$ does contain at
least two cycles.

\section{The group $A_{t,n}$}\label{sec6}
Similarly to \cite{rtv}, we  define a group $A_{t,n}$. Let $X=\{x,
y, z, \dots\}$ be a set of size $t$ and $R=\{r_x, r_y, \dots\}$ be
a set of size $t_1$, where $t_1\leq t$, and the indices of the
$r$'s are in a subset of $X$.

\begin{definition}
The group $A_{t,n}$ is generated by $(2n)^2|X|+2n|R|$ elements
$x_{i j}$, and $r_{x_k}$  where $x\in X$, \ \ $r\in R$ \ \ $i, j,
k \in \{1, 2, \dots, n, \bar{1}, \bar{2}, \dots, \bar{n} \}$ and
$\bar{\bar{i}}=i$ (we write $\bar{i}$ instead of $-i$).
\begin{eqnarray}
& x_{i i} =1 \label{26} \\
& x^{-1}_{i j}=x_{j i} \label{27} \\
& x_{i j}x_{j k}=x_{j k}x_{i j}= x_{i k} \ \ \mbox{for every $i$,
$j$ and $k$} \label{28} \\
& r_{x,i}r_{x,j}=x_{i \bar{j}} \ \ \mbox{for every $i$ and $j$} \label{29} \\
& x_{i j}y_{k l}= y_{k l}x_{i j} \ \ \ \mbox{and} \ \ \ x_{ij}x_{k
l}=x_{k l}x_{i j} \ \ \ \mbox{for every distinct $i, j, k, l$} \label{30} \\
\mbox{and in  addition} \nonumber \\
& x_{{\bar{j}} {\bar{i}}}=x_{i j} \label{31} \\
& x_{{\bar{i}} j}y_{{\bar{j}} k}x_{{\bar{k}} {\bar{i}}}y_{i
{\bar{j}}}x_{j{\bar{k}}}y_{k {\bar{i}}}=1 \label{32}  \\
& r_{x,i}y_{{\bar{i}} j}r_{x,j}r_{x,k}y_{{\bar{k}}
{\bar{i}}}r_{x,{\bar{i}}}r_{x,{\bar{j}}}y_{{j}
{\bar{k}}}r_{x,{\bar{k}}}=1 \label{33} \\
& r_{x,i}y_{{\bar{i}} j}r_{x,j}z_{{\bar{j}} k }r_{x,k}y_{{\bar{k}}
{\bar{i}}}r_{x,{\bar{i}}}z_{i {\bar{j}}}r_{x,{\bar{j}}}y_{j
{\bar{k}}} r_{x,{\bar{k}}}z_{k {i}}=1.\label{34}
\end{eqnarray}
\end{definition}

\begin{proposition}
For $n \geq 5$ or $t \leq 2$ the following (from \cite{rtv}) hold
in $A_{t,n}$:
\begin{eqnarray}
{}& [w_{i s},x_{j k}y_{k l}x_{k j}]=1  \ \ \mbox{for distinct $i,
j, k, l, s$} \label{pro1} \\
{}& x_{s i}y_{i j}x_{j s}w_{s k}=w_{s k}x_{k i}y_{i j}x_{j k} \ \
\mbox{for distinct $i, j, k, s$} \label{pro2}\\
{}& x_{s i}y_{i j}x_{j s}=x_{k i}y_{i j}x_{j k}  \ \ \mbox{for
distinct $i, j, k, s$} \label{ijks} \\
{}& [x_{s i}y_{i j}x_{j s},u_{j l}v_{l \bar{s}}u_{\bar{s} j}]=1 \
\ \mbox{for $t \leq 2$ or $n\geq 6$}. \label{pro3}
\end{eqnarray}
\end{proposition}

\begin{proof}
Relations (\ref{pro1}) and (\ref{pro2}) are proved in \cite{rtv}.

We prove Relation (\ref{ijks}). Let us consider the relation $x_{s
i}y_{i j}x_{j s}w_{s k}x_{k j} y_{j i}x_{i k}w_{k s}=1$. By
Relation (\ref{pro1}), this relation becomes $x_{s i}y_{i j}x_{j s}x_{k
j}y_{j i}x_{j k}w_{s k}x_{i j}w_{k s}=1$, and we are able to omit $w_{k
s}$ and $w_{s k}$ (since by Relation (\ref{30}), $w_{s k}x_{i
j}w_{k s}=x_{i j}w_{s k}w_{k s}=x_{i j}$). Therefore we get $x_{s
i}y_{i j}x_{j s}x_{k j}y_{j i}x_{j k}x_{i j}=1$, and this gives us
$(x_{s i}y_{i j}x_{j s})(x_{k j}y_{j i}x_{i k})=1$ (which is
exactly (\ref{ijks})).

Now we prove Relation (\ref{pro3}). If $n \geq 6$, there exist $t$
and $k$, distinct from $i, j, s, l$,  such that $x_{s i}y_{i
j}x_{j s}=x_{t i}y_{i j}x_{j t}$ and $u_{j l}v_{l
\bar{s}}u_{\bar{s} j}=u_{k l}v_{l \bar{s}}u_{\bar{s} k}$ (by
\ref{ijks}). And we can conclude that $[x_{t i}y_{i j}x_{j t},u_{k
l}v_{l \bar{s}}u_{\bar{s} k}]=1$ for $t \leq 2$ or $n\geq 6$.
\end{proof}

 \vskip 0.2cm

It is possible to define an action of $B_n$ on $A_{t,n}$ as
follows: $\sigma^{-1}x_{i j}\sigma:=x_{\sigma(i) \sigma(j)}$ and
$\sigma^{-1}r_{x,i}\sigma:=r_{x,\sigma(i)}$ for every $\sigma\in
B_n$ (similarly to the action of $S_n$ in \cite{rtv}).

\vskip 0.2cm

The $A_{t,n}$ has $t$ Abelian subgroups $Ab(x)$, where $Ab(x)$ is:
$x_{i j}$, $x_{\bar{i} j}$ for a particular $x$ or $x_{i j}$,
$x_{\bar{i} j}$ and $r_{x,k}$ for a particular $x$ (where
$r_{x,k}$ exists for the specific $x$ and $1\leq i,j, k\leq n$).
We see that the described groups $Ab(x)$ are abelian by using
Relations (\ref{28}), (\ref{29}), (\ref{30}) and (\ref{31}).

Each subgroup $Ab(x)$ is freely generated by $n$ elements $x_{i,
i+1}$ (where $1\leq i\leq n-1$) and $x_{\bar{1} 1}$ if $r_{x,j}$
does not exists. If $r_{x,j}$ exists, $Ab(x)$ is freely generated
by the $n$ elements $x_{i,i+1}$ (where $1\leq i\leq n-1$) and
$r_{x,1}$. In \cite[page 13]{rtv} it has been shown that the
subgroup $x_{i j}$, where $1\leq i,j\leq n$ is freely generated by
the set $x_{i,i+1}$, where $1\leq i\leq n-1$. Using Relation (\ref{28}),
$x_{\bar{i} j}=x_{\bar{i} \bar{1}}x_{\bar{1} 1}x_{1 j}$. Then using
Relation (\ref{31}), $_{\bar{i} j}=x_{1 i}x_{\bar{1} 1}x_{1 j}$. Hence,
adding a generator $x_{\bar{1} 1}$, we get all the elements $x_{i
j}\cup x_{\bar{i} j}$, where $1\leq i,j\leq n$. In case where
$r_{x,i}$ exists, by using Relation (\ref{29}), $r_{x_1}^2=x_{1
\bar{1}}$. Then using Relation (\ref{27}), $x_{1
\bar{1}}=x_{\bar{1} 1}^{-1}$. Hence, for $x$ where $r_{x,i}$
exists, $Ab(x)$ is freely generated by $x_{i,i+1}$ and $r_{x,1}$,
where $1\leq i\leq n-1$.

\section{The Main Theorem}\label{sec7}

\begin{theorem}\label{main}
Assume there is at least one anti-cycle or a loop in $T$. Then the
group $C_Y(T)$ is isomorphic to $A_{t,n} \rtimes D_n$ if there are no
loops in $T$. In the case of the existence of loops in $T$, it is
isomorphic to $A_{t,n} \rtimes B_n$.
\end{theorem}

In order to prove the theorem, we define, as in \cite{rtv},  a
spanning `tree' $T_0$. Note that for us, `tree' means that $T_0$
is connected and there are no cycles of any type in $T_0$ (no
cycles of $\tilde{A}$, $\tilde{B}$, $\tilde{C}$, $\tilde{D}$-type
).  But we allow the existence of anti-cycles (cycles with odd
number of edges, signed $-1$), and in particular we allow loops
and two conjugate edges to connect two vertices (which is an
anti-cycle of length $2$).

Now we explain how we get the spanning `tree' from the signed
graph of $C(T)$: In case of $\tilde{A}$-type, we get $T_0$ by
omitting one arbitrary edge, as it occurs in \cite{rtv}. In case
of $\tilde{D}$-type or $\tilde{B}$-type cycle, omitting one of the
edges in one of the anti-cycles (see Figures \ref{star1} and
\ref{star2}). In case of cycles of $\tilde{C}$-type, omitting one
of the loops $v$ or $w$ (see Figure \ref{dyn21}).

\vskip 0.2cm

We define $\gamma_i$ and $\gamma_{\bar{i}}$ for $1\leq i\leq n$. In case
of $\tilde{C}$-type cycle, where we omit a loop to get the spanning tree,
we define $\delta_{i}$ and $\delta_{\bar{i}}$ too.

We have already defined edges $u_i$ and $\bar{u_i}$, for every
$1\leq i\leq m$ in $\tilde{B}$, $\tilde{C}$, $\tilde{D}$-type
cycles with $m$ vertices. So, we can define certain elements
$\gamma_{i}$, $\gamma_{\bar{i}}$, $\delta_{i}$ and
$\delta_{\bar{i}}$ in every cycle in $T$:
\begin{eqnarray}
{}\gamma_{i}: &=& u_{i+2}u_{i+3}\cdots u_mu_1\cdots
u_i\\
{}\gamma_{\bar{i}}: & = &
u_{i+1}u_i\bar{u}_{i}u_{i+1}u_{i+2}\cdots
u_mu_1\cdots \bar{u}_{i-1}u_i\\
{}\delta_{i}: & = & u_{i}u_{i-1}\cdots u_{1}vu_{1}u_{2}\cdots
u_{m-1}wu_{m-1}u_{m-2}\cdots u_{i+1}\\
{}\delta_{\bar{i}}: & = & \delta_{i}^{-1}
\end{eqnarray}
for every $1\leq i\leq m$ and every $\tilde{C}$-type cycle of
length $m$ in $T$.

\vskip 0.2cm

Note that  the definition of $\gamma_{i}$ for $i>0$ is the same as
in \cite{rtv}.  In addition, we define $\gamma_{\bar{i}}$ too,
which has not been defined before. The following property is
important for the main theorem:
\begin{proposition}\label{7.1}
${\gamma_{\bar{i}}}^{-1}\gamma_{\bar{j}}={\gamma_{j}}^{-1}\gamma_{i}$
\ \ for every $i$ and $j$.
\end{proposition}

\begin{proposition}\label{7.2}
$[\gamma_{i}^{-1}\gamma_{j},\gamma_{k}^{-1}\gamma_{l}]=1$ for
every $i,j,k,l\in \{1, 2, \cdots m, \bar{1}, \bar{2}, \cdots,
\bar{m}\}$ \ \ ($i, j, k$ and $l$ are not necessarily distinct).
\end{proposition}

\begin{proposition}\label{7.3}
$[\delta_{i}, \gamma_{i}^{-1}\gamma_{j}]=1$ and
$\delta_{i}\delta_{j}=\gamma_{\bar{j}}^{-1}\gamma_{i}$.
\end{proposition}

\vskip 0.2cm \noindent {\bf{Proof of Propositions \ref{7.1},
\ref{7.2} and \ref{7.3}:}}
\\
The proof is by looking at the elements $\gamma_{i}$ (as it
defined) in the affine groups $\tilde{B}_m$, $\tilde{C}_m$ or
$\tilde{D}_m$ as periodic signed permutations with a period of
$2m+2$, which means $\pi(i+(2m+2))=\pi(i)+(2m+2)$ for every $i$
\cite{Eriksson}. Then for $j\neq \bar{i}$, the element
$\gamma_{j}^{-1}\gamma_{i}$ is the periodic signed permutation
$\pi$ which satisfies $\pi(i)=i+(2m+2)$, $\pi(j)=j-(2m+2)$,
$\gamma_{\bar{i}}^{-1}\gamma_{i}$ is the periodic signed
permutation $\pi(i)=i+2*(2m+2)$ and $\delta_{i}$ is the periodic
signed permutation in $\tilde{C}_m$ which satisfies
$\pi(i)=i+(2m+2)$. Since, $\pi(i)=i+(2m+2)$ means
$\pi(-i)=-i-(2m+2)$ and $pi(j)=j-(2m+2)$ means
$\pi(-j)=-j+(2m-2)$, Proposition \ref{7.1} holds. Propositions
\ref{7.2} and \ref{7.3} hold, since every two periodic
permutations $\pi$ and $\tau$ in an affine group $\tilde{B}_m$,
$\tilde{C}_m$ or $\tilde{D}_m$ commutes where $\pi(i)=i+k(2nm+2)$,
for every $i$ and some $k\in \mathbb{Z}$.

\vskip 0.4cm

\begin{proposition}
$T$ is a connected signed graph $T$ with $n$ vertices, then  it is
possible to extend the definition of $\gamma_{i}$ and
$\gamma_{\bar{i}}$ for every $1\leq i\leq n$.
\end{proposition}

\begin{proof}
The extension is done as follows. We define $\tilde{v}_i$ for
every edge $v_i$ in the signed graph $T$, in a similar way as it
was defined in \cite[page 7]{rtv}:

\vskip 0.2cm

$$
\tilde{v}_{a_v} = \left\{ \begin{array}{cc}
                   v_{a_v}, & \  \mbox{for every edge} \ \ v \
                   \mbox{signed by} \ \ a_v \ \mbox{which does not} \\
                    & \mbox{touch the cycle} \\
                   u_{i+1}, & \  \mbox{for every edge} \ \ v_{a_v}=u_i \\
                   \bar{u}_{i+1}, & \  \mbox{for every edge} \ \
                   v_{a_v}=\bar{u}_i \\
                   u_{i+1}v_{a_v}u_{i+1}, & \  \mbox{for every edge} \ \ v
                   \ \mbox{signed by} \ \  a_v \ \mbox{which does} \\
                    & \mbox{touch the cycle at vertex} \ \ x_i
                    \ \mbox{only} \\
                   u_{i+1}u_{j+1}v_{a_v}u_{j+1}u_{i+1}, & \  \mbox{for
                   every edge} \ \ v \ \mbox{signed by} \ \ a_v \
                   \mbox{which does} \\ & \mbox{touch the cycle at
                   vertices} \ \ x_i \ \mbox{and} \ \ x_j
                    \end{array} \right.
$$

$$\gamma_{t}:={\tilde{v}^{(s)}}_{a_{v^{(s)}}}\cdots
{\tilde{v}^{(1)}}_{a_{v^{(1)}}}\gamma_{1}{v^{(1)}}_{a_{v^{(1)}}}\cdots
{v^{(s)}}_{a_{v^{(s)}}}$$

and

$$\gamma_{\bar{t}}:={\tilde{v}^{(s)}}_{a_{v^{(s)}}}\cdots
{\tilde{v}^{(1)}}_{a_{v^{(1)}}}\gamma_{\bar{1}}{v^{(1)}}_{a_{v^{(1)}}}\cdots
{v^{(s)}}_{a_{v^{(s)}}}$$

where ${v^{(1)}}_{a_{v^{(1)}}}, \dots,
{v^{(s)}}_{a_{v^{(s)}}}$ is a connected path starting from the vertex $1$
and ending at the vertex $t$ in $T_0$.
\end{proof}

\begin{proposition}
We extend the definition of $\delta_{i}$ also to every $1\leq i\leq n$ in
case there is a loop $w\notin T_0$. $\delta_{i}$ has already been defined
for $\tilde{C}$-type cycle. Let $v^{(1)}\cdots v^{(s)}$ be a path from the
vertex $x_1$ in the $\tilde{C}$-type cycle to a vertex $x_t\in C_Y(T)$.
Then:

$$\delta_{t}:={v^{(s)}}_{a_{v^{(s)}}}\cdots
{v^{(1)}}_{a_{v^{(1)}}}\delta_{1}{v^{(1)}}_{a_{v^{(1)}}}\cdots
{v^{(s)}}_{a_{v^{(s)}}}$$

$$\delta_{\bar{t}}:={v^{(s)}}_{a_{v^{(s)}}}\cdots
{v^{(1)}}_{a_{v^{(1)}}}\delta_{\bar{1}}{v^{(1)}}_{a_{v^{(1)}}}\cdots
{v^{(s)}}_{a_{v^{(s)}}}$$

(Note: The definition of the extesion of $\delta$ is different from the
definition of the extension of $\gamma$, and does not use the
defined verices $\tilde{v}$).
\end{proposition}

\begin{remark}
We notice that Propositions \ref{7.1}, \ref{7.2} and \ref{7.3} are holds
for every $1\leq i\leq n$. The proof is by looking at the elements
$\gamma_{i}$ and $\delta_{i}$ as elements of the defined group, and
showing that the elements $\gamma_{i}^{-1}\gamma_{j}$ can be considered as
elements of the extended periodic permutation group to a period of $2n+2$,
where $\pi(j)=j+(2n+2)$ and $\pi(i)=i-(2n+2)$.
\end{remark}

\vskip 0.2cm
\noindent
{\bf Proof of Theorem \ref{main}:} \\
Similarly as defined in
\cite{rtv}, we define here $\theta : C_Y(T) \rightarrow A_{t,n}
\rtimes G$, where $t$ is the number of the cycles (every type) in
the signed graph, and $G= B_n$ or $D_n$, depending on existence of
loops in $T$.

For $u \in T$ we have $u = (ij)_a$ and
$$
\theta(u) =  \left\{ \begin{array}{cc}
            (ij)(\bar{i} \bar{j}), & \  \mbox{if} \ \  u \in T_0  \ \mbox{and} \ \ a=1 \\
            (i \bar{j})(\bar{i} j), & \  \mbox{if} \ \ u \in T_0 \
             \mbox{and} \ \ a= -1 \\
            (i j)(\bar{i} \bar{j}) v_{i j}, & \ \mbox{if}  \ \ v \notin
            T_0  \ \mbox{and} \ a=1 \\
            (i \bar{j})(\bar{i} j) v_{\bar{i} j}, & \ \mbox{if}  \ \ v \notin T_0 \ \mbox{and} \ a=-1 \ \
             \ \mbox{and} \\ & v \ \mbox{is not a loop in} \ \ \tilde{C}
             \ \mbox{type cycle} \\
             (i \bar{i})r_{v,i}, & \ \mbox{if}  \ \  v \notin T_0 \
             \mbox{and} \ \ v \ \mbox{is a loop} \\
             \end{array} \right.
$$

We can show that $\theta$ is well-defined on $C_Y(T)$, i.e., the
image of $\theta$ satisfies Relations (\ref{r1}), (\ref{r2}),
(\ref{r3})  and (\ref{r4}).
\begin{itemize}
\item
Relation (\ref{r1}) was treated in \cite{rtv}.
\item
(\ref{r2}) means that $\theta (uvu)$ commutes with $\theta (w
\bar{v} w)$ for every $u, v, w \in T$ and
$(uv)^3=(vw)^3=(uw)^2=1$. Now we treat the possible cases:
\begin{enumerate}
\item
$u, v \in T_0$: \ \ \ $\theta (uvu) =(i k)(\bar{i}\bar{k})$,
\item
$u \notin T_0$, $v \in T_0$: \\
$\theta (u) =(i j)(\bar{i} \bar{j})u_{i j}$, $\theta (v) =(k
j)(\bar{k} \bar{j})$, and $\theta (uvu) =(i j)(\bar{i}
\bar{j})u_{i j} (k j)(\bar{k} \bar{l})(i j)(\bar{i} \bar{j})u_{i
j} = (i k)(\bar{i} \bar{k})u_{i k}$,
\item
$u \in T_0$, $v \notin T_0$: \\
$\theta (u) =(i j)(\bar{i} \bar{j})$, $\theta (v) =(k j)(\bar{k}
\bar{j})v_{k j}$, and $\theta (uvu) =(i j)(\bar{i} \bar{j})(k
j)(\bar{k} \bar{j})v_{k j}(i j)(\bar{i} \bar{j})= (i k)(\bar{i}
\bar{k})v_{k i}$,
\item
$u, v \notin T_0$: \\
$\theta (u) =(i j)(\bar{i} \bar{j})u_{i j}$, $\theta (v) =(k
j)(\bar{k} \bar{j})v_{k j}$, and $\theta (uvu)=(i k)(\bar{i}
\bar{k})u_{j k}v_{k i}u_{i j}$ (see proof in \cite{rtv}).
\end{enumerate}
Similarly, $\theta (w \bar{v} w)$ is one of the followings:
\begin{enumerate}
\item $w, \bar{v} \in T_0$: \ \ \ $\theta (w \bar{v} w) =(l
\bar{j})(\bar{l} j)$,
\item $w \notin T_0$, $\bar{v} \in T_0$: \ \
\ $\theta (w \bar{v} w) =(l \bar{j})(\bar{l} j) w_{l \bar{j}}$,
\item $w \in T_0$, $\bar{v} \notin T_0$: \ \ \ $\theta (w \bar{v}
w) =(l \bar{j})(\bar{l} j) \bar{v}_{\bar{j} l}$, \item $w, \bar{v}
\notin T_0$: \ \ \ $\theta (w \bar{v} w) =(l \bar{j})(\bar{l} j)
w_{k \bar{j}} \bar{v}_{\bar{j} l} w_{l k}$.
\end{enumerate}

Since $i, k, l$ and $\bar{j}$ are distinct, each one of the elements  $(i
k)(\bar{i} \bar{k}), (i k)(\bar{i} \bar{k})u_{i k}, (i k)(\bar{i}
\bar{k})v_{k i}$ commutes with each one of the elements $(l
\bar{j})(\bar{l} j), (l \bar{j})(\bar{l} j)w_{l \bar{j}}, (l
\bar{j})(\bar{l} j)\bar{v}_{\bar{j} l}$.

It remains to show that each one of the elements $(i k)(\bar{i}
\bar{k}), (i k)(\bar{i} \bar{k})u_{i k}, (i k)(\bar{i}
\bar{k})v_{k i}$, and $(i k)(\bar{i} \bar{k})u_{j k}v_{k i}u_{i
j}$ commutes with  $(l \bar{j})(\bar{l} j)w_{k
\bar{j}}\bar{v}_{\bar{j} l}(j \bar{l})w_{l k}$. We start with (for
distinct $i, j, k, l$):
$$
(l \bar{j})(\bar{l} j)w_{k \bar{j}}\bar{v}_{\bar{j} l}w_{l k} (i
k)(\bar{i} \bar{k}) = (l \bar{j})(\bar{l} j)(i k)(\bar{i}
\bar{k})w_{i \bar{j}}\bar{v}_{\bar{j} l}w_{l i}= (i k)(\bar{i}
\bar{k})(l \bar{j})(\bar{l} j)w_{k \bar{j}}\bar{v}_{\bar{j} l}w_{l
k}, \mbox{by (\ref{pro2})}.
$$
Now we prove:
\begin{eqnarray*}
(l \bar{j})(\bar{l} j)w_{k \bar{j}}\bar{v}_{\bar{j} l}w_{l k} (i
k)(\bar{i} \bar{k})u_{i k} = (i k)(\bar{i} \bar{k})(l
\bar{j})(\bar{l} j)w_{i \bar{j}}\bar{v}_{\bar{j} l}w_{l i}u_{i
k}\\
\left . \begin{array}{c} = \\
[-.5cm]{\scriptstyle (\ref{pro2})}
\end{array}  \right .
(i k)(\bar{i} \bar{k})(l \bar{j})(\bar{l} j)u_{i {k}}w_{k
\bar{j}}\bar{v}_{\bar{j} l}w_{l k} = (i k)(\bar{i} \bar{k})u_{i
{k}}(l \bar{j})(\bar{l} j)w_{k \bar{j}}\bar{v}_{\bar{j} l}w_{l k}.
\end{eqnarray*}
Similarly,
$$
(l \bar{j})(\bar{l} j)w_{k \bar{j}}\bar{v}_{\bar{j} l}w_{l k} (i
k)(\bar{i} \bar{k})v_{k i} = (i k)(\bar{i} \bar{k})v_{k i}(l
\bar{j})(\bar{l} j)w_{k \bar{j}}\bar{v}_{\bar{j} l}w_{l k}.
$$

Now we show that if $n \geq 6$, then $(i k)(\bar{i} \bar{k})u_{j
k}{v}_{k i}u_{i j}$ commutes with $(l \bar{j})(\bar{l} j) w_{k
\bar{j}}\bar{v}_{\bar{j} l}w_{l k}$. Since $n \geq 6$, there exist
$p$ and $q$ such that $i, j, \bar{j}, k, l,p$ and $q$ are distinct
and by (\ref{ijks}), we have: $u_{j k}{v}_{k i}u_{i j} = u_{p
k}{v}_{k i}u_{i p}$ and $w_{k \bar{j}}\bar{v}_{\bar{j} l}w_{l k}=
w_{q \bar{j}}\bar{v}_{\bar{j} l}w_{l q}$.  Hence:
\begin{eqnarray*}
&(i k)(\bar{i} \bar{k})u_{j k}{v}_{k i}u_{i j}(l \bar{j})(\bar{l}
j)w_{k \bar{j}} \bar{v}_{\bar{j} l}w_{l k}
\left . \begin{array}{c} = \\
[-.5cm]{\scriptstyle (\ref{pro2})}
\end{array}  \right .
(i k)(\bar{i} \bar{k})u_{p k}{v}_{k i}u_{i p}(l \bar{j})(\bar{l}
j)w_{q \bar{j}} \bar{v}_{\bar{j} l}w_{l q} \\
&=(i k)(\bar{i} \bar{k})(l \bar{j})(\bar{l} j)u_{p k}{v}_{k i}u_{i
p}w_{q \bar{j}} \bar{v}_{\bar{j} l}w_{l q}
\left . \begin{array}{c} = \\
[-.5cm]{\scriptstyle (\ref{ijks})}
\end{array}  \right .
(i k)(\bar{i} \bar{k})(l \bar{j})(\bar{l} j)w_{q
\bar{j}}\bar{v}_{\bar{j} l}w_{l q}u_{p k}{v}_{k i}u_{i p} \\
&=(l \bar{j})(\bar{l} j)w_{q \bar{j}}\bar{v}_{\bar{j} l}w_{l q}(i
k)(\bar{i} \bar{k})u_{p k}{v}_{k i}u_{i p} = (l \bar{j})(\bar{l}
j)w_{q \bar{j}}\bar{v}_{\bar{j} l}w_{l q}(i k)(\bar{i}
\bar{k})u_{p k}{v}_{k i}u_{i p} \\
&\left . \begin{array}{c} = \\
[-.5cm]{\scriptstyle (\ref{pro2})}
\end{array}  \right . (l \bar{j})(\bar{l} j)w_{k \bar{j}}\bar{v}_{\bar{j} l}w_{l k}(i
k)(\bar{i} \bar{k})u_{j k}{v}_{k i}u_{i j}.
\end{eqnarray*}
\item
(\ref{r3}) means that $\theta (uvu)$ commutes with $\theta
(w)$ for $u, v, w \in T$ and $(uv)^3=(uw)^4=(vw)^4=1$, \ \ \
$\vcenter{\hbox{\epsfbox{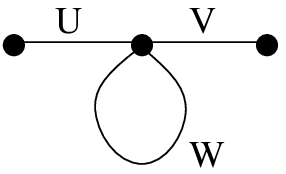}}}$.\\
The proof is the same one for the `fork' relation in \cite[P.
20]{rtv}. \item (\ref{r4}) means that $(\theta (u) \cdot
\theta(wvw))^3=1$ for $u, v, w \in T$ and
$(uv)^3=(uw)^4=(vw)^4=1$.
\end{itemize}

Now we classify the possible cases for $\theta(wvw)$ and $\theta(u)$. We
start with $\theta(wvw)$:
\begin{enumerate}
\item
$v, w \in T_0$: \\
$\theta (v) =(k j)(\bar{k} \bar{j})$, $\theta (w) =(\bar{j} j)$, and
$\theta (wvw) =(\bar{j} k)(j \bar{k})$.
\item
$v \notin T_0$, $w \in T_0$: \\
$\theta (v) =(k j)(\bar{k} \bar{j})v_{k j}$, $\theta (w) =(\bar{j}
j)$, and $\theta (wvw) =(\bar{j} k)(j \bar{k}) v_{k \bar{j}}$.
\item
$v \in T_0$, $w \notin T_0$: \\
$\theta (v) =(k j)(\bar{k} \bar{j})$, $\theta (w) =(\bar{j}
j)r_{w,\bar{j}}$, and $\theta (wvw) =(\bar{j} k)(j
\bar{k})r_{w,\bar{k}}r_{w,\bar{j}}$. \item
$v, w \notin T_0$: \\
$\theta (v) =(k j)(\bar{k} \bar{j})v_{k j}$, $\theta (w) =(\bar{j}
j)r_{w,\bar{j}}$, and $\theta (wvw) =(\bar{j} k)(j \bar{k})r_{w,\bar{k}}
v_{k \bar{j}}r_{w,\bar{j}}$.
\end{enumerate}
And here we give the following forms of $\theta({u})$:
\begin{enumerate}
\item[(a)]
$u \in T_0$: \ \ \ $\theta (u) =(i j)(\bar{i} \bar{j})$,\nonumber
\item[(b)]
$u \notin T_0$: \ \ \ $\theta (u) =(i j)(\bar{i} \bar{j})u_{i
j}$.\nonumber
\end{enumerate}

In the case (1) and (a) we have:
\begin{equation*}
(\theta (u) \cdot \theta(wvw))^3 = ((i j)(\bar{i} \bar{j})\cdot(\bar{j}
k)(j \bar{k}))^3=[(i \bar{k} j)(\bar{i} j \bar{k})]^3=1.
\end{equation*}
In the case (2) and (a) we have:
\begin{eqnarray*}
&&(\theta (u) \cdot \theta(wvw))^3 = (i j)(\bar{i} \bar{j})(\bar{j}
k)(j \bar{k})v_{k \bar{j}}(i j)(\bar{i} \bar{j})(\bar{j} k)(j \bar{k})v_{k
\bar{j}}(i j)(\bar{i} \bar{j})(\bar{j} k)(j \bar{k})v_{k \bar{j}}\\
&&=(i \bar{k} j)(\bar{i} j \bar{k})v_{k \bar{j}}(i \bar{k} j)(\bar{i}
k \bar{j})v_{k \bar{j}}(i \bar{k} j)(\bar{i} k \bar{j})v_{k
\bar{j}}=v_{\bar{i} k}v_{\bar{j} \bar{i}}v_{k \bar{j}}=1.
\end{eqnarray*}
In the case (3) and (a) we have:
\begin{eqnarray*}
&&(\theta (u) \cdot \theta(wvw))^3 =(i j)(\bar{i} \bar{j})(\bar{j}
k)(j \bar{k})r_{w,\bar{k}}r_{w,\bar{j}}(i j)(\bar{i} \bar{j})(\bar{j}
k)(j \bar{k})r_{w,\bar{k}}r_{w,\bar{j}}
(i j)(\bar{i} \bar{j})(\bar{j} k)(j \bar{k})r_{w,\bar{k}}r_{w,\bar{j}}\\
&&=(i \bar{k} j)(\bar{i} j \bar{k})r_{w,\bar{k}}r_{w,\bar{j}}(i
\bar{k} j)(\bar{i} j \bar{k})r_{w,\bar{k}}r_{w,\bar{j}}(i \bar{k}
j)(\bar{i} k
\bar{j})r_{w,\bar{k}}r_{w,\bar{j}}=r_{w,i}r_{w,k}r_{w,j}r_{w,\bar{i}}r_{w,\bar{k}}r_
{w,\bar{j}} =1.
\end{eqnarray*}
In the case (4) and (a) we have  (by Relation (\ref{32})):
\begin{eqnarray*}
(\theta (u) \cdot \theta(wvw))^3 =(i \bar{k} j)(\bar{i} j
\bar{k})r_{w,\bar{k}}v_{k \bar{j}}r_{w,\bar{j}} (i \bar{k}
j)(\bar{i} k \bar{j})r_{w,\bar{k}}v_{k \bar{j}}r_{w,\bar{j}} (i
\bar{k} j)(\bar{i} k \bar{j})r_{w,\bar{k}}v_{k \bar{j}}r_{w,\bar{j}}\\
=r_{w,i}v_{\bar{i} k}r_{w,k}r_{w,j}v_{\bar{j}\bar{i}}r_{w,\bar{i}}
r_{w,\bar{k}}v_{k \bar{j}}r_{w,\bar{j}}=1.
\end{eqnarray*}
In the case (1) and (b) we have (as in the case of (2) and (a)):
\begin{eqnarray*}
&&(\theta(wvw) \cdot \theta (u))^3 = ((k \bar{j})(\bar{k} j)(i
j)(\bar{i} \bar{j})u_{i {j}})^3=((i j \bar{k})(\bar{i} \bar{j} k)u_{i
{j}})^3=1.
\end{eqnarray*}
In the case (2) and (b) we have:
\begin{eqnarray*}
(\theta (u) \cdot \theta(wvw))^3 =((i j)(\bar{i} \bar{j})u_{{i} j}(j
\bar{k})(\bar{j} k)v_{k \bar{j}})^3=
((i \bar{k} j)(\bar{i} k \bar{j})u_{{i} \bar{k}}v_{k \bar{j}})^3=u_{j
{i}}v_{\bar{i} k}u_{\bar{k} j}v_{\bar{j} \bar{i}}
u_{i \bar{k}}v_{k \bar{j}}=1.
\end{eqnarray*}
In the case (3) and (b) we have:
\begin{eqnarray*}
&(\theta (u) \cdot \theta(wvw))^3 =((i j)(\bar{i} \bar{j})u_{{i}
j}(j \bar{k})(\bar{j} k)r_{w,\bar{k}}r_{w,\bar{j}})^3= ((i \bar{k}
j)(\bar{i} k \bar{j})u_{{i} \bar{k}}r_{w,\bar{k}}r_{w,\bar{j}})^3\\
&=u_{j {i}}r_{w,i}r_{w,k}u_{\bar{k} {j}} r_{w,j}r_{w,\bar{i}}u_{i
\bar{k}}r_{w,\bar{k}}r_{w,\bar{j}}=1.
\end{eqnarray*}
In the case (4) and (b) we have:
\begin{eqnarray*}
(\theta (u) \cdot \theta(wvw))^3 =((i j)(\bar{i} \bar{j})u_{{i}
j}(\bar{j} k)(j \bar{k})r_{w,\bar{k}}v_{k
\bar{j}}r_{w,\bar{j}})^3= ((i \bar{k} j)(\bar{i} k \bar{j})u_{{i}
\bar{k}}r_{w,\bar{k}}v_{k
\bar{j}}r_{w,\bar{j}})^3\\
=u_{j {i}}r_{w,i}v_{\bar{i} k} r_{w,k}u_{\bar{k}
j}r_{w,j}v_{\bar{j} \bar{i}}r_{w,\bar{i}}u_{i
\bar{k}}r_{w,\bar{k}} v_{{k} {j}}r_{w,\bar{j}}=1.
\end{eqnarray*}
We conclude that the Relations (\ref{r1}), (\ref{r2}), (\ref{r3})
and  (\ref{r4}) are satisfied for $\theta(C_Y(T))$. Hence $\theta$
is well defined on $C_Y(T)$.

\vskip 0.2cm

This proves $\theta : C_Y(T) \rightarrow A_{t,n} \rtimes G$ is
a homomorphism.  $G=Im(\phi(C(T))$, where  $\phi$ is the natural map
from $C(T)$ into $B_n$ which was defined in Lemma \ref{natural}. By the
same Lemma $G= B_n$ in case $T$ does contain a loop, or $G=D_n$ in
case $T$ does not contain a loop but does conatain an anti-cycle.

\vskip 0.2cm

Now we define
$\tau : A_{t,n}
\rtimes G \rightarrow C_Y(T)$ as it was defined in \cite{rtv}:

$\tau(v)=v$ if $v\in B_n$ or $v\in D_n$,
$\tau(x_{i j})={\gamma_j}^{-1}{\gamma_i}$ and
$\tau(r_{x,i})=\delta_{i}$.

We need to check Relations (\ref{26}), (\ref{27}), (\ref{28}) and
(\ref{30}) for $\tau(x_{i,j})$. Relation (\ref{26}) holds
trivially since $\gamma_{i}^{-1}\gamma_{i}=1$. Relation (\ref{27})
holds, since
$(\gamma_{i}^{-1}\gamma_{j})^{-1}=\gamma_{j}^{-1}\gamma_{i}$, and
Relation (\ref{28}) and (\ref{30}) hold by Proposition \ref{7.2}.

Hence, $\tau$ is well defined, and $\tau$ is the inverse map of
$\theta$. There are five options:

1) $\theta \tau(x_{i j})=\theta(\gamma_{j}^{-1}\gamma_{i})$. Then the
proof is exactly the same as in \cite{rtv}.

2) $\theta \tau(x_{\bar{i}
\bar{j}})=\theta(\gamma_{\bar{j}}^{-1}\gamma_{\bar{i}})$. By Proposition
\ref{7.1}: \
$\gamma_{\bar{j}}^{-1}\gamma_{\bar{i}}=\gamma_{i}^{-1}\gamma_{j}$,
and $\tau(\gamma_{\bar{j}}^{-1}\gamma_{\bar{i}})=x_{i
j}=x_{\bar{j} \bar{i}}$, by Relation~(\ref{31}).

3) $\theta \tau(x_{\bar{i} j})=\theta(\gamma_{j}^{-1}\gamma_{\bar{i}})=$\\
$=\theta(u_{j}u_{j-1}\cdots u_{1}u_{m}x_{0 m}u_{m-1}\cdots
u_{j+2}u_{i+1}u_{i}\bar{u}_{i}u_{i+1}\cdots u_{m}x_{0 m}u_1\cdots
u_{i-2}\bar{u}_{i-1}u_{i})=$\\
$=x_{\bar{i} j}$.

4) $\theta \tau(x_{i \bar{j}})=\theta(\gamma_{\bar{j}}^{-1}\gamma_{i})=
\theta(\gamma_{i}^{-1}\gamma_{\bar{j}})^{-1}=x_{\bar{j}
i}^{-1}=x_{i\bar{j}}$ by Relation (\ref{27}).

5) $\theta \tau({r_{x,i}})= \theta(\delta_{i})=
\theta(u_{i}u_{i-1}\cdots u_{1}vu_{1}\cdots u_{m-1}wu_{m-1}\cdots
u_{i+1})=$ \\
$=u_{i}u_{i-1}\cdots u_{1}vu_{1}\cdots
u_{m-1}wr_{x,{m-1}}u_{m-1}\cdots u_{i+1}=r_{x,i}$.

Hence, in every case $\theta \tau$ is the identity, then $\tau$ is
the inverse map to $\theta$. $\square$

\end{document}